\documentclass{amsart}
\usepackage[raiselinks,colorlinks]{hyperref}
\usepackage{a4,enumerate}
\usepackage{amssymb}
\usepackage{exscale}
\newcommand{\be}{\begin{equation}}
\newcommand{\ee}{\end{equation}}
\newcommand{\bea}{\begin{eqnarray*}}
\newcommand{\eea}{\end{eqnarray*}}
\newcommand{\beq}{\begin{eqnarray}}
\newcommand{\eeq}{\end{eqnarray}}
\renewcommand{\r}{\right}
\renewcommand{\l}{\left}

\newcommand{\cC}{\mathcal C}
\newcommand{\cD}{\mathcal D}

\newcommand{\cF}{\mathcal F}
\newcommand{\cG}{\mathcal G}
\newcommand{\cH}{\mathcal H}
\newcommand{\cI}{\mathcal I}
\newcommand{\CC}{\mathbb{C}}
\newcommand{\EE}{\mathbb{E}}
\newcommand{\NN}{\mathbb{N}}
\newcommand{\PP}{\mathbb{P}}
\newcommand{\QQ}{\mathbb{Q}}
\newcommand{\RR}{\mathbb{R}}
\newcommand{\ZZ}{\mathbb{Z}}
\newtheorem{thm}{Theorem}

\newtheorem{lem}[thm]{Lemma}
\newtheorem{cor}[thm]{Corollary}
\newtheorem{hyp}{Assumption}
\theoremstyle{remark}
\newtheorem{rem}[thm]{Remark}
\newtheorem{dfn}[thm]{Definition}

 \newcommand{\G}{G}
 \newcommand{\E}{E}
 
 \newcommand{\tr}{\text{tr\,}}

\newcommand{\ph}{{\varphi}}

\newcommand{\Nach}{{\,\rightarrow\,}}

\newcommand{\sk}{{\,|\,}}

\newcommand{\Hm}[1]{\leavevmode{\marginpar{\tiny%
$\hbox to 0mm{\hspace*{-0.5mm}$\leftarrow$\hss}%
\vcenter{\vrule depth 0.1mm height 0.1mm width \the\marginparwidth}%
\hbox to
0mm{\hss$\rightarrow$\hspace*{-0.5mm}}$\\\relax\raggedright #1}}}

\begin{document}

\title[Wegner estimate for Schr\"odinger operators on metric graphs]
{A linear Wegner estimate for alloy type Schr\"odinger operators on metric graphs}

\author[M.~Helm]{Mario Helm$^{1}$}

\address[$^1$]{Fakult\"at f\"ur Mathematik,\, 09107\, TU-Chemnitz, Germany   }
\email{ Mario.Helm@mathematik.tu-chemnitz.de}

\author[I.~Veseli\'c]{Ivan Veseli\'c$^{1,2}$}
\address[$^2$]{Emmy-Noether-Programme of the Deutsche Forschungsgemeinschaft}
\urladdr{\href{http://www.tu-chemnitz.de/mathematik/schroedinger/members.php}{www.tu-chemnitz.de/mathematik/schroedinger/members.php}}



\begin{abstract}
We study spectra of alloy-type random Schr\"odinger operators on
metric graphs. For finite edge subsets of general graphs we prove
a Wegner estimate which is linear in the volume (i.e.~the number
of edges) and the length of the considered energy interval. The
single site potential of the alloy-type model needs to have fixed
sign, but the considered metric graph does not need to have a
periodic structure. The second result we obtain is an exhaustion
construction of the integrated density of states for ergodic
random Schr\"odinger operators on metric graphs with a
$\ZZ^{\nu}$-structure. For certain models the two above results
together imply the Lipschitz continuity of the integrated density
of states.
\end{abstract}

\maketitle

\section{Introduction}

In the present paper we study spectral properties of random Schr\"odinger operators
on a metric graph. More precisely, we consider a countable metric graph and a random
Hamiltonian with a  so-called alloy type potential. Under suitable assumptions we are
able to prove a Wegner estimate for the restriction of the random Hamiltonian to a
finite part of the graph. Our Wegner estimate is optimal in the sense that it
is linear in the energy and the volume.

In the case that the integrated density of states (IDS) of the
random Schr\"odinger operator exists, the Wegner estimate implies
that the IDS is Lipschitz continuous. For certain random operators
on $\ZZ^{\nu}$-metric graphs we establish the existence of the IDS via
an exhaustion procedure. More generally, the existence of an
selfaveraging IDS is ensured by certain amenability assumptions on
the metric graph and ergodicity assumptions on the random
Hamiltonian, see Remark \ref{r-genIDS}.

The literature on Wegner estimates for Schr\"odinger operators in
the continuum, i.e. on $\RR^{\nu}$, is quite extensive. We refer
to the textbook accounts
\cite{CyconFKS-87,CarmonaL-90,PasturF-92,Stollmann-01,Veselic-04a,KirschM}
and the references therein. For the application to spectral
localization in random media, a quite weak form of Wegner estimate
is sufficient. However, upper bounds which are linear both in the
volume and the energy interval length are of independent interest,
since they imply the regularity of the IDS. To obtain such a
Wegner estimate, the proofs which are presented in the literature
need some assumptions, beyond those necessary for a weak form of
Wegner estimate. In particular, mostly the covering condition
\[
\sum_{k \in \ZZ^{\nu}} u(x-k) \geq const. > 0 \quad \text{ for all } x\in
\RR^{\nu}
\]
was assumed. Here $u$ is a compactly supported, non-negative
single site potential. Papers which are devoted to the question
how this assumption can be removed or at least weakened include
\cite{Klopp-95a,Kirsch-96,Veselic-96,CombesHN-01,KirschV-02b,Giere-02,CombesHK-03}.
In the recent \cite{CombesHK} a linear Wegner estimate is
proposed without the abovementioned covering condition. However,
the deterministic part of the Schr\"odinger operator needs to be
translation invariant. The reason for this assumption is that it
implies a uniform, translation invariant kind of unique
continuation principle for the solutions of the Schr\"odinger
equation, cf.~Section 4 in \cite{CombesHK-03}.

\medskip

In our setting, on metric graphs, the situation is somewhat simpler, since
in certain aspects the Schr\"odinger operators behave as in the one dimensional case.
In particular, on the (one dimensional) edges one has a stronger (and simpler) unique continuation principle at disposal.
We use techniques developed for Wegner estimates for operators on $L^2(\RR)$ in
\cite{KirschV-02b}. (See also \cite{Veselic-96,Giere-02,CombesHK-03} for similar results
for one dimensional random Schr\"odinger operators.)

This allows us to prove linear Wegner estimates for random
operators on very general metric graphs. In particular, we do not
need any periodicity condition on the graph. Our assumptions are general enough to include, for instance,
operators on the metric graph associated to $\ZZ^{\nu}$ or, more
generally, the Cayley graph of a finitely generated
discrete group, or a Penrose tiling graph.

Note that for quantum graphs there exists a certain dichotomy concerning the unique continuation principle.
While along a single edge one has a quite strong version of this property, globally on the graph it does not hold.
More precisely, a Laplacian (even without potential) on a metric graph with ``free boundary conditions'' at the vertices
may exhibit compactly supported eigenfunctions. For ergodic Hamiltonians on discrete graphs it is well known (see for instance
\cite{DodziukLMSY-03,KlassertLS-03,Veselic-05b,Veselic-06,LenzV}) that compactly supported eigenfunctions are related to discontinuities of the IDS.
Thus, for a complete understanding of the continuity properties of the IDS a systematic analysis of compactly supported
eigenfunctions is necessary. See also \cite{GruberLV} for related results and discussion.
\medskip

The structure of the paper is as follows: in the next section we formulate our results,
Section \ref{s-proofWegner} contains the bulk of the proof of the Wegner estimate,
while arguments relating to boundary condition perturbations are deferred to Section \ref{s-bc},
and those to the unique continuation principle to Section \ref{s-smallu}. The last section contains a
proof of the construction of the IDS via an exhaustion procedure.

\section{Model and Results}
\label{s-results}

First we define the geometric structure of the metric graphs and the operators acting on the associated $L^2$-Hilbert space.

\begin{dfn}
Let $V$ and $E$ be countable sets and $\cG$ a map
\[
\cG \colon E \to V\times V \times [0,\infty), \quad e \mapsto
(\iota(e),\tau(e),l_e).
\]
We call the triple $G=(V,E,\cG)$ a metric graph, elements of $V=V(\G)$ vertices, elements of $E=E(\G)$ edges,
$\iota(e)$ the initial vertex of $e$, $\tau(e)$ the terminal vertex of $e$ and
$l_e$ the length of $e$. Both $\iota(e)$ and $\tau(e)$ are called endvertices of $e$,
or incident to $e$. The two endvertices of an edge are allowed to coincide.
The number of edged incident to the vertex $v$ is called the
\textit{degree of $v$}. We assume that the degree is finite for all vertices.
The elements of $V$ which have degree equal to one we call \textit{boundary vertices}
and denote their set by $V^\partial$. The elements of $V^i:=V\setminus V^\partial$
are called \textit{interior vertices}.

We will consider the graph $\G=(V,E, \cG)$ as a topological space, more
precisely as an one dimensional CW-complex, which we will again
denote by $\G$. The $0$-skeleton of $\G$ is $V$ and the collection
of its one dimensional cells is given by $E$. Each one-dimensional
cell $e \in E$ is attached either to one or to two zero dimensional cells
$v \in V$, namely $\iota(e)$ and $\tau(e)$. This defines the
topological structure of $G$.
\end{dfn}
\smallskip

We identify each edge $e$ with the open interval $\,(0,l_e)$,
where the point $0$ corresponds to the vertex $\iota(e)$ and $l_e$
to $\tau(e)$. Assume that there exist constants $0<l_-,
l_+<\infty$ such that for all $e\in \E$
\begin{equation*}
    0<l_- \leq l_e \leq l_+ <\infty .
\end{equation*}
The identification of edges by intervals allows us to define in a
natural way the length of a path between two points in the
topological space $G$. Taking the infimum over the lengths of
paths connecting two given points in $G$, one obtains a distance
function $d \colon G\times G\to [0,\infty)$. Since we assumed that
each vertex of $G$ has bounded degree, the map $d$ is indeed a
metric, cf.~for instance Section 2.2 in \cite{Schubert-06}. Thus
we have turned $G$ into a metric space $(G,d)$.

\smallskip

For a finite subset $\Lambda \subset E$ we define the subgraph
$\G_\Lambda $ by deleting all edges $e \in E\setminus \Lambda$
and the arising isolated vertices. We denote the set of vertices of
$G_\Lambda$ by $V_\Lambda$, the set of vertices $v\in V_\Lambda$ of degree one (in $G_\Lambda$)
by $V_\Lambda^\partial$, and its complement $V_\Lambda \setminus V_\Lambda^\partial$
by $V_\Lambda^i$. Again, elements of $V_\Lambda^\partial$
are called boundary vertices of $G_\Lambda$
and elements of  $V_\Lambda^i$ interior vertices of $G_\Lambda$.
Similarly as above we may consider $\G_\Lambda$ as a sub-CW-complex of $\G$ with an induced topology and metric.

\smallskip

For any $\Lambda \subset E$ the Hilbert spaces $L_2(\G_\Lambda)$
have a natural direct sum representation
$L_2(\G_\Lambda)=\oplus_{e\in \Lambda}L_2(0,l_e)$. In particular
for $\Lambda = E$ we have $L_2(\G)= \oplus_{e\in \E}L_2(0,l_e)$,
and for $\tilde\Lambda\subset \Lambda \subset E$ we have
$L_2(\G_{\tilde\Lambda})=\oplus_{e\in
{\tilde\Lambda}}L_2(0,l_e)\subset L_2(\G_\Lambda)=\oplus_{e\in
\Lambda}L_2(0,l_e)$.

For a function $\phi\colon \G \to \CC$ and an edge $e\in E$ we denote by
$\phi_e:=\phi|_e$ its restriction to $e$ (which is identified with $(0,l_e)$).
We denote by $C(\G)$ the space of continuous, complex-valued functions on the metric space $(\G,d)$.
Similarly, $C(\G_\Lambda)$ denotes the space of continuous, complex-valued functions on the metric sub-space $(\G_\Lambda,d)$.
For each $v \in V$, any edge $e$ incident to $v$, and function $f\in W^{2,2}(e) \subset C^1(e)\cong C^1(0,l_e)$
we define the derivatives
\begin{align}
\partial_e f(v) &:= \partial_e f(0)   :=  \lim_{\epsilon\searrow 0} \frac{f(\epsilon)-f(0)}{\epsilon} & \text{ if } v = \iota(e)
\intertext{and}
\partial_e f(v) &:= \partial_e f(l_e) :=  \lim_{\epsilon\searrow 0} \frac{f(l_e)-f(l_e-\epsilon)}{\epsilon} & \text{ if } v = \tau(e).
\end{align}
Note that, since $f|_e\in W^{2,2}(e)$, by the Sobolev imbedding
theorem the function $f$ is not only continuously differentiable
on the open segment $(0,l_e)$, but also that its derivative has
well defined limits at both boundaries $0$ and $l_e$.

For any $\Lambda \subset E$ it will be convenient to use the
following Sobolev space
\begin{align*}
W^{2,2}(\Lambda):= \oplus_{e\in \Lambda } W^{2,2}(e) \subset
C^1(\Lambda) := \oplus_{e\in \Lambda }  C^{1}(e)
\\
\text{ with the norm }\|\phi\|^2_{W^{2,2}(\Lambda)}:= \sum_{e\in \Lambda} \|\phi_e\|^2_{W^{2,2}(0,l_e)} .
\end{align*}

Note that this space is defined on the edge set only and does not
see the graph structure of $G$. The operators which we consider will be defined on the space
\begin{align} \nonumber
\cD(\Delta_\Lambda):=
\{f & \in W^{2,2}(\Lambda) \cap C(G_\Lambda)
\mid
\\
&\forall v \in V_\Lambda^i :\sum_{e\in E, \iota(e)=v}\!\!\!\!
\partial_e f(v)=\sum_{e\in E, \tau(e)=v} \!\!\!\!\partial_e f(v), \quad
\forall v \in V_\Lambda^\partial : f(v) =0\}.
\end{align}
\smallskip

For each $\Lambda \subset E$ we define a linear operator
\[
-\Delta_\Lambda \colon \cD(\Delta_\Lambda) \to L^2(\Lambda)
\]
by the rule
\[
(-\Delta_\Lambda f)(x) := -\frac{\partial^2 f_e(x)}{\partial x^2 }
\]
if $x \in G$ is contained in the edge $e$. This way the function
$-\Delta_\Lambda f$ is defined on the set $E\subset G$, whose complement $V=G\setminus E$
in the metric space $G$ has Hausdorff measure zero.

The operator $-\Delta_\Lambda$ is selfadjoint on the domain
$\cD(\Delta_\Lambda)$, see for instance \cite{Kuchment-04,KostrykinS-99b} or Section 3.3 in \cite{Schubert-06}.
At the boundary vertices it has clearly Dirichlet boundary conditions,
while the type of boundary conditions it has at the interior vertices is called
``free boundary conditions'' by some authors and ``Kirchhoff boundary conditions'' by others.

\smallskip

Next we describe the potential energy part of the Hamiltonian.
Since it is random, we need an appropriate probability space.

Let $\omega_- < \omega_+$ be real numbers, $\Omega$ a probability space,
and $\PP$ a probability measure on $\Omega$.
Denote by $\EE$ the mathematical expectation on $\Omega$ with respect to $\PP$.
Let $W\colon \Omega \times G \to [\omega_-,\omega_+]$ be
a stochastic process which is jointly measurable in the variables $\omega\in \Omega $
and $x\in G$. For a fixed $\omega \in \Omega$ we denote
by $W(\omega)\colon L_2(G) \to L_2(G)$ the multiplication operator
$W(\omega)f(x) =(W(\omega,x)f)(x)$, $x \in G$.

To derive our Theorems \ref{t-Wegner} and \ref{t-IDS} below we will need more specific
hypotheses on  $\PP$ and $W(\omega)$. These are formulated in what follows.
The first assumption describes random potentials of alloy type for which we are
able to prove a Wegner estimate.

\begin{hyp}
\label{a-Wegner} Let $c_+ \ge c_->0$, $s>0$ and $c_g$ be  real
numbers. For each edge $e\in \E$ let $\mu_e$ be a probability
measure which is absolutely continuous with respect to Lebesgue
measure. More precisely, for $e\in \E$ let $g_e\in
L_\infty[\omega_-, \omega_+]$, $\|g_e\|_\infty \leq c_g$, and
$d\mu_e(t)=g_e(t)\,dt$ on the interval $[\omega_-, \omega_+]$. Let
the probability space $\Omega$ be given by the cartesian product
$\times_{e \in \E}[\omega_-,\omega_+]$, and $\PP$ by the product
$\otimes_{e\in \E} \mu_e$.

For each $e\in \E$ let $u_e\in L_\infty(0,l_e)$ be a \emph{single site potential} satisfying
\begin{equation*}
    c_-\, \chi_{S_e}  \leq u_e \leq c_+\, \chi_{[0,l_e]},
\end{equation*}
where $S_e\subset [0,l_e]$ is an interval of length $ |S_e|\geq s$.
We imbed $L_\infty(0,l_e)\cong L_\infty(e)$ in $L_\infty (\G)$
and thus consider $u_e$ as an element of the latter space.
Let the random potential $W(\omega)$ have the form
\begin{equation*}
    W(\omega)\colon L_2(G) \to L_2(G) , \quad W(\omega)= \sum_{e\in \E} \omega_e
    u_e.
\end{equation*}
\end{hyp}

For quantum graphs with a $\ZZ^{\nu}$-structure we will establish the existence
of the integrated density of states in Theorem \ref{t-IDS}.
The following hypothesis formulates precisely what type of $\ZZ^{\nu}$-structure we need for this result.

\begin{hyp}
\label{a-IDS} The vertex set $V$ consists of the points $\ZZ^{\nu}
\subset \RR^{\nu}$. The set of edges $E$ consists of segments
$e=[x,y]$ parallel to the coordinate axes in $\RR^{\nu}$ where $x,y\in
V$ have Euclidean distance equal to one. The union $V \cup E \subset \RR^{\nu}$ inherits
the metric structure of $\RR^{\nu}$. This structure coincides with the one defined earlier
and we denote the corresponding metric space again by $(G,d)$.

Assume that for each $x\in \ZZ^{\nu}$  there is a measure preserving map $\tau_x\colon \Omega\to\Omega$ such that
$\{\tau_x\}_{x\in\ZZ^{\nu}}$ forms an additive group which acts ergodically on $(\Omega, \cF,\PP)$.

Assume that the random potential transforms in the following way under translations:
$W(\tau_x\omega, y)=W(\omega, y-x)$ for all $x\in \ZZ^{\nu}, y \in G$ and $\omega \in \Omega$.
Thus the Schr\"odinger operator is \emph{equivariant} in the sense that
\begin{equation}
\label{e-equivariant}
U_x H(\omega)U_x^* = H(\tau_x\omega) \quad \text{ for all $x\in \ZZ^{\nu}$ and $\omega \in \Omega$}
\end{equation}
where $U_x f(y)=f(y-x), y \in G$ denotes the unitary translation operator by $x \in \ZZ^{\nu}$.
\end{hyp}

The restriction of the operator  $W(\omega)$ to $L_2(\Lambda)$
will be denoted by $W_\Lambda(\omega)$ or, if the set $\Lambda$ is
clear from the context, simply by $W(\omega)$ again. Note that the
norm of the operator $W_\Lambda(\omega)$  is uniformly bounded in
$\omega \in \Omega$ and $\Lambda \subset E$. Thus the operator
$H_\Lambda(\omega):= -\Delta_\Lambda+ W_\Lambda(\omega)$ is
selfadjoint on $\cD(\Delta_\Lambda)$ for any $\omega \in \Omega$
and $\Lambda \subset E$. If the set $\Lambda\subset E$ contains
a single element $e$ we write $-\Delta_e(\omega)$ for
$-\Delta_\Lambda(\omega)$ and $H_e(\omega)$ for
$H_\Lambda(\omega)$.
\bigskip

A metric graph $G$ together with a Schr\"odinger $H$ operator which is defined on it
we call a \textit{quantum graph}.

\bigskip

Our main result is the following:
\begin{thm}
\label{t-Wegner}
If Assumption \ref{a-Wegner} holds,
then there exists for any $\lambda\in \RR$  a constant $C$ such that for all
$\epsilon\in [0,1]$,
\begin{equation*}
  \EE\{ \tr \chi_{[\lambda-\epsilon,\lambda+\epsilon]}(H_\Lambda(\omega)) \} \leq C \cdot\epsilon\cdot \sharp\Lambda .
\end{equation*}
\end{thm}

If in the terminology of Assumption \ref{a-Wegner} we have $S_e=
[0,l_e]$ for all $e\in E$, then the constant appearing in Theorem
\ref{t-Wegner} can be chosen uniformly in the energy parameter
$\lambda\in \RR$.
\bigskip

Now we turn to the situation where the quantum graph has a $\ZZ^{\nu}$-structure
as formulated in Assumption \ref{a-IDS}, and describe how the
integrated density of states can be defined by an exhaustion procedure.

For any $l\in \NN$ denote by $\Lambda_l$ the set of edges which are contained in
\[
\{x \in \RR^{\nu} \mid x_i \in (0,l) \text{ for all } i =1,\dots,d
\},
\]
and abbreviate $V_{\Lambda_l}$ by $V_l$, and $G_{\Lambda_l}$ by
$G_l$. Denote by $\partial G_l$ the boundary of $G_l$ as a subset
of the metric space $G$, i.e.~the vertices in $V_l$ which have
degree one. Similarly abbreviate $\Delta_l := \Delta_{\Lambda_l}$
and $ H_l(\omega):=H_{\Lambda_l}(\omega)$. These are restrictions
of the operators $\Delta$ and $H(\omega)$ to a finite cube with
sidelength $l$ and with Dirichlet boundary conditions on the
boundary of the cube. The finite cube Schr\"odinger operator
$H_l(\omega)$ is again self-adjoint, lower bounded and has purely
discrete spectrum. Let us enumerate the eigenvalues of
$H_l(\omega)$ in ascending order $\lambda_1(H_l(\omega)) <
\lambda_2(H_l(\omega)) \le \lambda_3(H_l(\omega)) \le \dots$ \,
and counting multiplicities.
\medskip

Thus for each $\lambda \in \RR$ and $l \in \NN$, the counting function
\[
F_\omega^l(\lambda) := \sharp \{n \in \NN \mid \lambda_n(H_l(\omega)) \le \lambda\}
\]
is monotone increasing and right-continuous, i.e.~a distribution function,
which is associated to a pure point measure. Denote by $N_\omega^l(\lambda) := \frac{1}{l^{\nu}} F_\omega^l(\lambda)$
the volume-scaled version of $F_\omega^l(\lambda)$.
\medskip


\begin{thm}
\label{t-IDS}
Let Assumption \ref{a-IDS} hold, then there exists a distribution function $N\colon \RR \to \RR$
and a subset $\Omega' \subset \Omega$ of measure one such that for all $\omega \in \Omega$ and for all
$\lambda\in \RR$ where $N$ is continuous the convergence
\begin{equation}
\label{e-convergence}
\lim_{l \to \infty} N_\omega^l(\lambda) = N(\lambda)
\end{equation}
holds.

\end{thm}

\begin{rem}
\label{r-genIDS} Under certain additional assumptions it is
possible to prove that the normalized finite volume eigenvalue
counting functions converge uniformly in the energy parameter to
the IDS, see \cite{GruberLV}. For metric graphs whose isometry
group does not exhibit a $\ZZ^{\nu}$-structure, but is merely
amenable, it should be still possible to define the IDS of an
ergodic Hamiltonian by an exhaustion, proceeding analogously as in
the paper \cite{LenzPV-04}.
\end{rem}

\section{Proof of the Wegner estimate}
\label{s-proofWegner}

For the purposes of the proof it will be necessary to differentiate the spectral projection
with respect to the energy parameter, which motivates the introduction of the following
smooth 'switch function'.

Let  $\rho$ be a smooth, non-decreasing function such that on
$(-\infty, -\epsilon]$ it is identically equal to $-1$, on
$[\epsilon, \infty)$ it is identically equal to zero and
$\|\rho'\|_\infty \le 1/\epsilon$. Then
\[
\chi_{(\lambda-\epsilon,\lambda+\epsilon )} (x) \le
\rho(x-\lambda+2\epsilon) -\rho(x-\lambda-2\epsilon) =
\int_{-2\epsilon}^{2\epsilon}  \,\rho'(x-\lambda+t) \,dt
\]
Thus by the spectral theorem
\[
\chi_{(\lambda-\epsilon,\lambda+\epsilon )}(H_\Lambda(\omega)) \le
\int_{-2\epsilon}^{2\epsilon} dt
\,\rho'(H_\Lambda(\omega)-\lambda+t)
\]
in the sense of quadratic forms. Since $B_\epsilon(\lambda)
=(\lambda-\epsilon,\lambda+\epsilon )$ is bounded and $
\sigma(H_\Lambda(\omega))$ discrete, the above operators are trace
class and we may estimate:
\[
\tr \Big [ \chi_{B_\epsilon(\lambda)}(H_\Lambda(\omega)) \Big ]
\le \tr \Big [ \int_{-2\epsilon}^{2\epsilon}
\,\rho'(H_\Lambda(\omega)-\lambda+t) \,dt \Big ] = \sum_{n\in \NN}
\int_{-2\epsilon}^{2\epsilon}
\,\rho'(\lambda_n^\Lambda(\omega)-\lambda+t)\,dt
\]
where $\lambda_n^\Lambda(\omega)$ denotes the eigenvalues of $H_\Lambda(\omega)$ enumerated in
non-decreasing order and counting multiplicities. Only a finite number of terms in the sum are non-zero.

\smallskip

\begin{proof}[Proof of Theorem \ref{t-Wegner}]
Let $\rho$ be as above.  Similarly as in \cite{Kirsch-96}, p. 509, we
estimate
\begin{equation*}
    \EE \{ \tr \chi_{[\lambda-\epsilon,\lambda+\epsilon]}(H_\Lambda(\omega)) \} \leq \int\limits_{[\omega_-,
    \omega_+]^{\Lambda}}  \!\sum_{n\in
    \NN}\, \int\limits_{-2\epsilon}^{2\epsilon} \rho'(\lambda_n^\Lambda(\omega)-\lambda+ t)\,dt
    \prod_{\tilde{e}\in \Lambda}
    d\mu_{\tilde{e}}.
\end{equation*}
To bound the right hand side we use Corollary \ref{f999} in Section \ref{s-smallu} and the monotone
convergence theorem to obtain the estimate
\begin{eqnarray}\label{e-997}
    \lefteqn{C_1 \int\limits_{[\omega_-, \omega_+]^{\Lambda}}  \sum_{n\in\NN}\,
    \int\limits_{-2\epsilon}^{2\epsilon}  \sum_{{e}\in \Lambda}
    \frac{\partial \rho(\lambda_n^\Lambda(\omega)-\lambda+ t)}{\partial\omega_{{e}}}\,dt \prod_{\tilde{e}\in \Lambda} d\mu_{\tilde{e}} (\omega_{\tilde e})}
     \nonumber\\
    &&{}\leq\, C_1\sum_{e\in \Lambda}   \int_{[\omega_-, \omega_+]^{A}}
    \prod_{\tilde{e}\in A}
    \!d\mu_{\tilde{e}}(\mu_{\tilde e})  \sum_{n\in\NN}\, \int\limits_{-2\epsilon}^{2\epsilon}\int\limits_{\omega_-}^{\omega_+}
    \frac{\partial \rho(\lambda_n^\Lambda(\omega)-\!\lambda\!+ t)}{\partial\omega_e}\,d\mu_{e}(\omega_{e})dt.
\end{eqnarray}
Here we used the abbreviation $A :=\Lambda \setminus \{e \}$.
Denote by $H_\Lambda(e,\omega_+)$ and $H_\Lambda(e,\omega_-)$ the
operator $H_\Lambda$, where the random variable $\omega_{e}$ is
set to its maximum respectively its minimum value. The $n-$th eigenvalues
of this operators are abbreviated by
$\lambda_n^\Lambda(e,\omega_+)$ and
$\lambda_n^\Lambda(e,\omega_+)$. Using monotonicity, the sum over
$n$ in (\ref{e-997}) can be estimated by
\begin{equation}\label{e-996}
   \|g_{e}\|_\infty \sum_{n\in \NN}\, \int\limits_{-2\epsilon}^{2\epsilon}
   \{\,\rho (\lambda_n^\Lambda(e,\omega_+)-\lambda+t)- \rho (\lambda_n^\Lambda(e,\omega_-)-\lambda+t)\, \}\,dt.
\end{equation}
We introduce now the operators $H_\Lambda^\ast(e,\omega),\,
\ast\in \{D,N\} $, that coincide with $H_\Lambda(\omega)$ up to
additional Dirichlet, respectively Neumann b.c. at the vertices
$\iota(e)$ and $\tau(e)$. Their eigenvalues are denoted by
$\lambda_n^{\Lambda,\ast}$.
By Lemma \ref{l-DNbracketing} on Dirichlet-Neumann bracketing
in Section \ref{s-bc} we have
\begin{multline}
 \label{e-995}
 \rho  (\lambda_n^\Lambda(e,\omega_+)-\lambda+t)- \rho  (\lambda_n^\Lambda(e,\omega_-)-\lambda+t)
 \\ \le
 \rho  (\lambda_n^{\Lambda,D}(e,\omega_+)-\lambda+t)- \rho  (\lambda_n^{\Lambda,N}(e,\omega_-)-\lambda+t).
\end{multline}
Because of the decoupling of the edge $e$, the latter operators
can be written as the direct sum $H_\Lambda^\ast(e,\omega)=
H_\Lambda^{e,\ast}(e,\omega) \oplus H_\Lambda^{c,\ast}(e,\omega)$
of operators acting on $L_2(0,l_{e})$ and $L_2(\G_\Lambda
\!\setminus\! e)$, respectively. 
For the eigenvalues of these operators we use the notation 
$\lambda_n^{e,*}(e,\omega)$ and $\lambda_n^{c,*}(e,\omega)$.
Hence the sum over the terms in (\ref{e-995})
can be separated in the corresponding parts:
\begin{eqnarray}
 \lefteqn{\sum_{n\in \NN} \rho  (\lambda_n^{\Lambda,D}(e,\omega_+)-\lambda+t)-
 \rho  (\lambda_n^{\Lambda,N}(e,\omega_-)-\lambda+t)}
 \\ \label{e-edge-ev}
  &=&  \sum_{n\in \NN}  \rho  (\lambda_n^{e,D}(e,\omega_+)-\lambda+t)-
  \rho (\lambda_n^{e,N}(e,\omega_-)-\lambda+t)
 \\ \label{e-993}
 &+&  \sum_{n\in \NN} \rho (\lambda_n^{c,D}(e,\omega_+)-\lambda+t)-
 \rho (\lambda_n^{c,N}(e,\omega_-)-\lambda+t).
\end{eqnarray}
We estimate first the sum in \eqref{e-edge-ev}. The difference
$H_\Lambda^{e,D}- H_\Lambda^{e,N}$ is a perturbation of rank 2
in resolvent sense, see for instance \cite{Simon-95}.
By Lemma  \ref{l-finite rank} in Section \ref{s-bc}
the first term in \eqref{e-edge-ev} obeys the bound
\begin{equation*}
\rho  (\lambda_n^{e,D}(e,\omega_+)-\lambda+t)
\le
\rho  (\lambda_{n+2}^{e,N}(e,\omega_+)-\lambda+t).
\end{equation*}
By a telescoping argument we obtain the estimate
\bea
  \sum_{n\in \NN}  \rho  (\lambda_{n+2}^{e,N}(e,\omega_+)-\lambda+t)- \rho
 (\lambda_n^{e,N}(e,\omega_+)-\lambda+t) \le 2
\eea
where we used that the total variation of $\rho$ equals one.
Thus we are left with estimating
\begin{multline}
\label{e-traces}
  \sum_{n\in \NN}   \rho (\lambda_n^{e,N}(e,\omega_+)-\lambda+t) -
  \rho (\lambda_n^{e,N}(e,\omega_-)-\lambda+t)
\\
=
  \tr \big [\rho (H_\Lambda^{e,N}(e,\omega_+)-\lambda+t) -
  \rho (H_\Lambda^{e,N}(e,\omega_-)-\lambda+t) \big].
\end{multline}
Note that by definition of $\rho$ we have
$\rho(x-\lambda +t)\le \chi_{[\lambda-3\epsilon,\infty)}(x) $ and
$-\rho(x-\lambda +t)\le \chi_{(\lambda+3\epsilon,\infty)}(x) $ for all $t \in [-2\epsilon,2\epsilon]$.
Thus the trace in \eqref{e-traces} is bounded by
\begin{eqnarray*}
  \lefteqn  {\tr [\chi_{[\lambda-3\epsilon,\infty)}(H_\Lambda^{e,N}(e,\omega_+)) -  \chi_{(\lambda+3\epsilon,\infty)}
  (H_\Lambda^{e,N}(e,\omega_-))]  }\\
  & \leq&  \tr [\chi_{[\lambda-3\epsilon,\infty)}(H_\Lambda^{e,N}(e,\omega_-)+
  (\omega_+-\omega_-)\|u_{e}\|_\infty) -
  \chi_{(\lambda+3\epsilon,\infty)}(H_\Lambda^{e,N}(e,\omega_-))] \\
  &=&  2 + \tr (\chi_{[\lambda-3\epsilon - (\omega_+-\omega_-)\|u_{e}\|_\infty ,
  \lambda+3\epsilon)}(H_\Lambda^{e,N}(e,\omega_-)))\\
  &\le&   C_2 .
\end{eqnarray*}
Here $C_2$ is a constant which is independent of $\Lambda$ and
depends only on the considered energy interval. In fact, a careful
look reveals that $C_2$ depends only on the length $(6\epsilon +
(\omega_+-\omega_-)\|u_{e}\|_\infty)$ of the energy interval, but
not on $\lambda \in \RR$. See for instance Section 3 in \cite{GruberLV}.

Next we want to estimate the sum in \eqref{e-993}.
Let $\tilde{d}:= \deg(\iota(e))+ \deg(\tau(e))-2$. Similarly as above,
one sees that the difference $H_\Lambda^{c,D}(e,\omega_+)-H_\Lambda^{c,N}(e,\omega_+)$
is a  perturbation of rank $\tilde{d}$ in resolvent sense.
Consequently, the first term in \eqref{e-993} obeys the bound
\begin{equation*}
   \rho (\lambda_{n}^{c,D}(e,\omega_+)-\lambda+t)\leq \rho
 (\lambda_{n+\tilde{d}}^{c,N}(e,\omega_+)-\lambda+t).
\end{equation*}
A telescoping argument bounds the whole sum in (\ref{e-993}) by $\tilde{d}$ times
the total variation of $\rho$, i.e. by $\tilde{d}$.

Hence the following upper bound on \eqref{e-997} completes the proof:
\begin{equation*}
   C_1 \,\!\! \sum_{e\in \Lambda}\!\|g_{e}\|_\infty \, \int_{-2\epsilon}^{2\epsilon}
   (\tilde{d}+C_2)\, dt
   \,\leq\, 4 C_1 (C_2+\tilde{d})\, c_g \cdot \epsilon \cdot \sharp \Lambda.
\end{equation*}

\end{proof}

\section{Changing boundary conditions}
\label{s-bc}

In this section we show how to control the shifting of eigenvalues
when one modifies boundary conditions. This is used in the
proof of Theorem \ref{t-Wegner}  above.

Lemma \ref{l-finite rank} shows that  changing the boundary
conditions from the  Dirichlet operator
$H_\Lambda^{e,D}(e,\omega_-)$ to the Neumann one
$H_\Lambda^{e,N}(e,\omega_-)$ shifts the eigenvalue index at most
by the rank of the perturbation, i.e. by $2$. The proof of this
fact is based on the min-max formula for eigenvalues, see
e.g.~\cite{ReedS-78}, Theorem XIII.1.

Lemma \ref{l-DNbracketing} is a monotonicity statement about
boundary conditions known as Dirichlet-Neumann-bracketing. By
quadratic form considerations one sees that the passage from
Neumann to Dirichlet boundary condition shifts the eigenvalues up.

First we proof an auxiliary lemma for bounded operators.

\begin{lem}\label{l-034}
Let $S,T=S+V$ be bounded, selfadjoint operators on a separable Hilbert
space $\cH$ with dim $(\text{ker }V)^\perp = d < \infty$. Assume
that the spectra of $S$ and $T$ below $\inf
\sigma_{\text{ess}}(S)$, respectively $\inf
\sigma_{\text{ess}}(T)$, consist of an infinite, discrete set of
eigenvalues.

Then for the $m$-th eigenvalues of $S$ and $T$, counted in ascending order including multiplicities,
we have $\lambda_m(T)\leq \lambda_{m+d}(S)\quad (\forall m\in \NN)$.
\end{lem}
\begin{proof} Let $\tilde{L}:= (\text{ker }V)^\perp$. By the min-max
principle one has
\begin{eqnarray*}
  \lambda_m(T) &=& \max_{\dim L= m-1} \,\,\min_{\|\ph\|=1,\, \ph\perp L} \,(\ph\sk T\ph) \\
   &\leq&  \max_{\dim L= m-1} \,\,\min_{\|\ph\|=1,\, \ph\perp
   L+\tilde{L}}\,
   [(\ph\sk S\ph)+\underbrace{(\ph\sk V\ph)}_{=0}]\\
   &\leq&  \max_{\dim \hat{L}= m+d-1} \,\,\min_{\|\ph\|=1,\, \ph\perp \hat{L}}\,(\ph\sk S\ph) \\
   &=& \lambda_{m+d}\,(S).
\end{eqnarray*}
\end{proof}

The lemma above can of course not be applied to the Hamiltonians under
consideration because they are unbounded operators. But there are some
related operators for which the statement is true - one can simply
compare suitable resolvents and rearrange the eigenvalues by the
spectral theorem.

\begin{lem} \label{l-finite rank}
For the operators $H_\Lambda^{e,D}(e,\omega_-)$ and
$H_\Lambda^{e,N}(e,\omega_-)$ we have
\begin{equation}\label{g988}
   \lambda_m(H_\Lambda^{e,D}(e,\omega_-))\leq
   \lambda_{m+2}(H_\Lambda^{e,N}(e,\omega_-)).
\end{equation}
\end{lem}

\begin{proof}
We shift the spectrum of both operators by addition of a suitable
constant and work in the following with the two arising strictly
positive operators $H_1$ (Dirichlet case) and $H_2$ (Neumann
case).

Let $D_0:= D(H_1)\cap D(H_2)$. Then one has
\begin{equation}\label{e-990}
    H_1^{-1} - H_2^{-1} \,|_{H_2 D_0} = H_1^{-1} (H_2-H_1)
    H_2^{-1} \,|_{H_2 D_0}.
\end{equation}
By definition, $H_1|_{D_0}-H_2|_{D_0}=0$, such that by (\ref{e-990}) and
continuity we get
\begin{equation*}
       H_1^{-1} - H_2^{-1} \,|_{\overline{H_2 D_0}} = 0,
\end{equation*}
i.e. ker\,$(H_1^{-1} - H_2^{-1})^\perp \subset H_2 D_0^\perp$.

We want to apply Lemma \ref{l-034} to the bounded operators
$-H_1^{-1}$ and $-H_2^{-1}$. So we have to show that $\dim \ker
(H_1^{-1}-H_2^{-1})^\perp \leq 2$, for what in turn $\dim H_2
D_0^\perp \leq 2$ is sufficient.
For this purpose we show next that $D_0$ and $D(H_2)$ differ  by a
$2$-dimensional subspace $L$.

Let $\phi_1,\phi_2 \in C^\infty(0,l_e)$ be such that $\phi_1\equiv
1$ in a neighborhood of $0$ and $\phi_1\equiv 0$ in a neighborhood
of $l_e$, and similarly let $\phi_2\equiv 0$ in a neighborhood of
$0$ and $\phi_2\equiv 1$ in a neighborhood of $l_e$. Thus
$\phi_1,\phi_2$ are linearly independent vectors and moreover
elements of $D(H_2)\setminus D_0$.

Let $\psi$ an arbitrary element of $D(H_2)$, $c_1:= \psi(0)$ and
$c_2:= \psi(l_e)$. Then
\begin{equation*}
  (\psi-c_1\phi_1-c_2\phi_2)'(0) \quad\text{ and } \quad(\psi-c_1\phi_1-c_2\phi_2)'(l_e)=0 \\
\end{equation*}
as well as
\begin{equation*}
   (\psi-c_1\phi_1-c_2\phi_2)(0)=0 \quad \text{ and } \quad(\psi-c_1\phi_1-c_2\phi_2)(l_e)= 0.
\end{equation*}
Hence $\psi- c_1\phi_1 - c_2\phi_2 \in D_0$ what implies that $\dim
D(H_2)\setminus D_0 = 2$.

Now, $H_2: D(H_2)\Nach \mathcal{H}$ is one to one, and we get
\begin{eqnarray*}
   \mathcal{H} &=& H_2 D(H_2) \\
   &=&  H_2 (D_0+ L) \\
   &=&  H_2 D_0 + H_2 L.
\end{eqnarray*}
Hence $ \dim H_2 D_0 ^\perp \leq \dim H_2 L \leq 2$.

So Lemma \ref{l-034} is applicable to $-H_1^{-1}$ and $-H_2^{-1}$,
and inequality (\ref{g988}) follows by the spectral theorem.
\end{proof}
It can  be seen easily that the operators
$H_\Lambda^{c,D}(e,\omega_-)$ and $H_\Lambda^{c,N}(e,\omega_-)$ in
the proof of Theorem \ref{t-Wegner} can be treated in the same
way.
\smallskip

\begin{lem}
\label{l-DNbracketing} Let $\omega \in \Omega$, $\Lambda \subset
E$ finite, and $e \in \Lambda$ be arbitrary. Consider the
operators $H_\Lambda(e,\omega)$, $H_\Lambda^D(e,\omega)$  and
$H_\Lambda^N(e,\omega)$ and their eigenvalues
$\lambda_n^{\Lambda}(e,\omega)$, $\lambda_n^{\Lambda,D}(e,\omega)$
and $\lambda_n^{\Lambda,N}(e,\omega)$ (defined in Section
\ref{s-proofIDS}). Then the following inequalities hold  for all $n
\in \NN$ \beq \lambda_n^{\Lambda,N}(e,\omega) \le
\lambda_n^\Lambda(e,\omega)  \le \lambda_n^{\Lambda,D}(e,\omega).
\eeq
\end{lem}
\begin{proof}
Denote by $h_\Lambda(e,\omega)$, $h_\Lambda^D(e,\omega)$ and $h_\Lambda^N(e,\omega)$ the  quadratic forms
associated to the operators $H_\Lambda(e,\omega)$, $H_\Lambda^D(e,\omega)$  and $H_\Lambda^N(e,\omega)$.
Then by Section 4.2 in \cite{Schubert-06} the quadratic form domains obey
\[
h_\Lambda^N(e,\omega) \supset h_\Lambda(e,\omega) \supset h_\Lambda^D(e,\omega).
\]
Thus the statement of the lemma follows immediately if one applies the
quadratic form version of the min-max formula for eigenvalues, see e.g.~\cite{ReedS-78},
Theorem XIII.2.
\end{proof}

\section{Single site potentials of small support}\label{s-smallu}

In this section we prove a uniform lower bound on the sum of derivatives
of eigenvalues which is formulated in
\begin{lem}
\label{l-lowerbound}
\begin{equation}
\label{Eder} \sum_{e \in \Lambda} \frac{\partial
\lambda_n^\Lambda(\omega)}{\partial \omega_e} \ge C_1(I) >0
\end{equation}
for all eigenvalues $\lambda_n^\Lambda$ of $H_\Lambda(\omega)$
inside a bounded energy interval $I$. The bound $C(I)$ does not
depend on the set of edges  $\Lambda \subset E$ and on the eigenvalue index
$n\in \NN$.
\end{lem}

We infer immediately:

\begin{cor} \label{f999}
Let $\rho: \RR\to [0,1]$ be a smooth, monotone function with
$\rho=-1$ on $(-\infty, -\epsilon]$ and $\rho=0$ on $[\epsilon,
\infty)$. Then, for the $n-$th eigenvalue $\lambda_n^\Lambda$ of
$H_\Lambda(\omega)$ we have
\begin{equation*}
    \rho'(\lambda_n^\Lambda(\omega)-\lambda+ t)\, \leq \, C_1(I)
    \sum_{e\in \Lambda} \frac{\partial \rho(\lambda_n^\Lambda(\omega)-\lambda+ t)}{\partial
    \omega_e}.
\end{equation*}
\end{cor}
\begin{proof}

By the chain rule
\begin{eqnarray*}
 \sum_{ e \in \Lambda}
\frac{ \partial \rho (\lambda_n^\Lambda(\omega) -\lambda + \theta) }{ \partial \omega_e }
=
\rho'(\lambda_n^\Lambda(\omega) -\lambda + \theta)
\sum_{e \in \Lambda} \frac{\partial \lambda_n^\Lambda(\omega)}{\partial \omega_e} \, ,
\end{eqnarray*}
\eqref{Eder} implies the estimate:
\begin{eqnarray*}
\rho'(\lambda_n^\Lambda(\omega) -\lambda + \theta) \le C_1(I)^{-1}
 \sum_{ e \in \Lambda}
\frac{ \partial \rho (\lambda_n^\Lambda(\omega) -\lambda + \theta)
}{ \partial \omega_e } .
\end{eqnarray*}
\end{proof}

To infer the lower bound \eqref{Eder}, set $S= \bigcup_{e \in
\Lambda} \ S_e$ and apply the Hellman-Feynman theorem, i.e.~first
order perturbation theory. For a normalized eigenfunction $\psi_n$
corresponding to $\lambda_n^\Lambda(\omega)$ we have:
\begin{eqnarray*}
\sum_{e \in \Lambda} \frac{\partial
\lambda_n^\Lambda(\omega)}{\partial \omega_e} = \sum_{e \in
\Lambda} (\psi_n \sk u_e\psi_n) \ge \int_{S} |\psi_n|^2.
\end{eqnarray*}
If the integral on the rightern side would extend over the whole
of $\Lambda$, it would be equal to $1$ due to the normalization of
$\psi_n$. A priori the integral over $S$ could be much smaller,
but the following Lemma shows that we can control the ratio of the
two integrals.
\begin{lem}
\label{uc} Let $I$ be a bounded interval and $S_e \subset [0,l_e]$
a non-degenerate interval. There exists a constant $C_1(I)>0$ such
that \be \label{e-uc} \int_{S_e} |\psi|^2 \ge  C_1(I)\int_0^{l_e}
|\psi|^2 \ee for  all $\Lambda\subset E$ finite, all $e \in
\Lambda$ and for any eigenfunction $\psi$ corresponding to an
eigenvalue $\lambda\in I$ of $H_\Lambda(\omega)$.
\end{lem}
Thus $\int_{S} |\psi|^2 \ge  C_1(I)\int_{\G_\Lambda} |\psi|^2$ with the same constant as in \eqref{e-uc}.
Hence Lemma \ref{uc} implies directly Lemma \ref{l-lowerbound}.
\begin{proof}
For $y\in \RR$ denote by $S_e+y:=\{x \in (0,l_e)\mid x-y \in S_e\}$
the translations of the set $S_e$ along the edge $e$. The derivative of the function
\begin{equation}
 \phi ( y) := \int_{S_e+y}  |\psi(x)|^2 \,dx = \int_{S_e}
 |\psi(x-y)|^2 \,dx
\end{equation}
satisfies
\begin{eqnarray*}
\l | \frac{\partial}{\partial y}  \phi ( y)   \r | & = & \l |
\int_{S_e}  \l [  \frac{\partial}{\partial y} \psi (x-y)  \r ]
\overline{ \psi(x-y) } \,dx
        + \int_{S_e}  \psi(x-y) \, \frac{\partial}{\partial y} \overline{ \psi(x-y) } \,dx  \r |
\\
& \leq & 2 \l \|  \psi \r \|_{ L^2( S_e+y ) } \l \| \psi'\r \|_{ L^2( S_e+y ) }.
\end{eqnarray*}
Sobolev norm estimates (e.g.~Theorems 7.25 and 7.27 in \cite{GilbargT-83}) imply
\[
\|  \psi' \|_{L^2( S_e+y )} \le  C_3 \,  \|  \psi \|_{L^2( S_e+y )} +\|  \psi'' \|_{L^2( S_e+y )}.
\]
By the eigenvalue equation  we have
\begin{equation}
\l | \frac{\partial}{\partial y}  \phi ( y)   \r |
\le C_4 \  \|  \psi  \|_{ L^2( S_e+y ) }^2 = C_4 \, \phi ( y), \qquad C_4 =C_4(\|W-\lambda\|_\infty).
\end{equation}
Gronwall's Lemma  implies  $ \phi(y)  \le  \exp(C_4 |y| ) \, \phi(0)$ and thus
\begin{eqnarray*}
\int_0^{l_e} |\psi|^2 \le  e^{C_4 l_e}    \ \frac{l_e}{|S_e|} \int_{S_e} |\psi|^2.
\end{eqnarray*}
\end{proof}

\section{Proof of Theorem \ref{t-IDS}}
\label{s-proofIDS}

The proof of Theorem \ref{t-IDS} on the existence of the
integrated density of states follows the arguments of
\cite{KirschM-82c}. There the convergence \eqref{e-convergence} of
the IDS was established for all rational energies $\lambda \in
\QQ$,  whereas we establish the same fact for all $\lambda$ which
are continuity points of the IDS. The main step consists in
proving that a superadditive ergodic theorem  of the paper
\cite{AkcogluK-81} is applicable.

In the whole of this section we assume that the quantum graph has a $\ZZ^{\nu}$-structure as described in Assumption \ref{a-IDS}.

\medskip

Let us first describe the type of superadditive processes
considered in \cite{AkcogluK-81}. Denote by $T$ the semigroup of
measurable transformations given by $\tau_x,\, x \in \NN_0^{\nu}$,
where $\NN_0:= \NN\cup \{0\}$, and by $\cI$ the class of sets of
the form
\[
\{x \in \RR^{\nu}\mid a_i \le x_i < b_i, \text{ for all } i= 1, \dots, d\}
\]
where $a, b \in \NN_0^{\nu}$. For $x\in \NN_0^{\nu}$ and ${Q} \in \cI$
denote by ${Q}+x:=\{y \in \RR^{\nu} \mid y-x \in {Q}\}$ the translation of the set ${Q}$ by $x$.

\begin{dfn}
\label{d-superadditive}
A function $F\colon \cI \to L^1(\Omega, \PP)$ which satisfies
\begin{enumerate}[(i)]
\item $F_{Q} \circ \tau_x = F_{{Q}+x}$ for all ${Q} \in \cI$ and $x \in \NN_0^{\nu}$,
\item if ${Q}_1, \dots, {Q}_n\in \cI$ are disjoint sets and
if their union ${Q}=\cup_{i=1}^n {Q}_i$ is again an element of $\cI$, then
\[
F_{Q} \ge \sum_{i=1}^n F_{{Q}_i},
\]
\item $\gamma(F):= \sup\limits_{{Q} \in \cI, |{Q}|>0} \frac{1}{|{Q}|} \, \EE\{F_{Q}\}< \infty$
\end{enumerate}
is called a (discrete) \emph{superadditive process}.

\end{dfn}

We state Theorem 2.4 from \cite{AkcogluK-81}:
\begin{thm}
\label{t-AkcogluKrengel} If $F$ is a discrete superadditive
process and ${Q}_l:= \{x \mid 0 \le x_i < l\linebreak \text{ for
all } i= 1, \dots, d\}$, then
\[
\lim_{l \to \infty} \frac{F_{{Q}_l}(\omega)}{l^{\nu}}
\]
exists for almost all $\omega \in \Omega$.
\end{thm}
In fact, in the case that the semigroup $T$ acts ergodically on the probability space $(\Omega, \PP)$
one can identify the limit, see the remark on page 59 in  \cite{AkcogluK-81}:
\[
\lim_{l \to \infty}   l^{-d}\, F_{{Q}_l}(\omega) = \gamma(F) \quad \text{ almost surely.}
\]
In Section 6.2 of \cite{Krengel-85} the above statements are
extended to the case that it is not the semigroup $T$ which acts
ergodically on $\Omega$, but rather the full group $\tau_x,\, x
\in \ZZ^{\nu}$.
\bigskip

To apply the superadditive ergodic theorem we consider for arbitrary,
fixed $\lambda\in \RR$ the process given by the eigenvalue counting functions
of the Schr\"odinger operator with Dirichlet 
boundary conditions.
For $Q \in \cI$ denote by $\Lambda$ the set of edges $e \in E \subset \RR^{\nu}$
which are contained in the \emph{interior} of $Q$ and set
\begin{align*}
F_Q:= F_{Q}(\lambda,\omega)&:= \sharp\{n| \,
\lambda_n(H_\Lambda(\omega)) < \lambda\}, \quad Q \in \cI.
\end{align*}
The Dirichlet Schr\"odinger operator $H_\Lambda(\omega)$ has been defined in Section \ref{s-results}.
We will show that $F_Q,\, Q \in \cI,$ is a superadditive process.

\begin{lem}
\label{l-DNbWa}
For any fixed energy value $\lambda$, the process $F_Q$ is superadditive.

\end{lem}
\begin{proof}
We have to check that the properties (i) -- (iii) in Definition \ref{d-superadditive} hold.
\smallskip

Since $\Lambda = \text{int\,} Q \cap E$ we have also $\Lambda +x=
(\text{int\,}  Q +x) \cap E$ for all $x \in \ZZ^{\nu}$. The
equivariance property \eqref{e-equivariant} of the random
operators carries over to the spectral projections and thus to the
eigenvalue counting functions. Hence property (i) holds.

\smallskip

Property (ii) can bee seen using Lemma \ref{l-DNbracketing} on  Dirichlet-Neumann bracketing in Section \ref{s-bc}.
For $Q$ and $ Q_1, \dots, Q_n
\in \cI$ as in (ii), the set $\Lambda$ contains $\bigcup_{i=1}^n
\Lambda_i$ and finitely many edges $e_1, \dots, e_N$ which lie in
$\text{int\,} Q \setminus \bigcup_{i=1}^n \text{int\,} Q_i$. By
introducing Dirichlet boundary conditions at finitely many
vertices one obtains from the operator $H_\Lambda(\omega)$ the
direct sum
\begin{equation}
\label{e-dirsum}
\bigoplus_{i=1}^n H_{\Lambda_i}(\omega)  \, \oplus \, \bigoplus_{j=1}^N H_{e_j}(\omega).
\end{equation}
Hence the eigenvalue counting function of $H_\Lambda(\omega)$
is an upper bound of the one of the direct sum operator \eqref{e-dirsum}. Obviously the counting function of
$\bigoplus_{i=1}^n H_{\Lambda_i}(\omega)$
can be estimated from above by the one of the operator \eqref{e-dirsum}. Now property (ii) follows.
\smallskip

Denote the eigenvalue counting function of the negative Dirichlet
Laplacian on $e=(0,1)$ by
\[
n_0(\lambda) := \sharp \{n \in \NN \mid \lambda_n(-\Delta_{e}) \le \lambda\},
\]
Obviously the counting function of
\begin{equation}
\label{e-dirsumedges} -\bigoplus_{e \in \Lambda} \Delta_e
\end{equation}
equals $|\Lambda| \cdot n_0(\lambda)$. Since the operators
$-\Delta_\Lambda$ and \eqref{e-dirsumedges} differ by Dirichlet
boundary conditions at $2 \, |\Lambda|$, or even less, vertices, the
counting function of $-\Delta_\Lambda$ is bounded by
\[
|\Lambda| \cdot n_0(\lambda)  + 4 d \cdot |\Lambda|.
\]
For this conclusion we use an argument analogous to  Lemma
\ref{l-finite rank} in Section \ref{s-bc}.

Since the random potentials we are considering are uniformly bounded by a constant, say $K$, we have
\[
F_Q(\lambda, \omega) \le |\Lambda| \cdot n_0(\lambda+K)  + 4 d
\cdot |\Lambda| \quad \text{ for all } \lambda \in \RR,\, Q \in
\cI.
\]
The number of edges in the set $\Lambda$ which is associated to a box $Q \in \cI$
is linearly bounded by the volume of $Q$.
Hence (iii) is proven.

\end{proof}
\bigskip
Now we can complete the
\begin{proof}[Proof of Theorem \ref{t-IDS}]
For a fixed $\lambda \in \RR$ one applies the ergodic theorem of
\cite{AkcogluK-81} to the superadditive process
$F_Q(\lambda,\omega),\, Q \in \cI$. Let us denote the
corresponding $\gamma(F)$ by $\gamma(\lambda)$. By definition
$F_Q(\lambda,\omega)\le F_Q(\tilde\lambda,\omega)$ for all
$\lambda \le \tilde \lambda \in \RR$ and all $\omega\in \Omega$,
$Q \in \cI$, thus $\lambda \mapsto \gamma(\lambda)$ is a
non-decreasing function. In particular, it has at most a countable
set of discontinuity points. Denote the complement of this set by
$\cC$ and choose a dense countable set $S\subset \cC$. Hence
$\gamma$ is continuous at each $\lambda \in S$.

Since in our case the transformation group is ergodic,
for each $\lambda$ there is a set $\Omega_\lambda$ of measure one on which the convergence
$\lim_{l \to \infty} l^{-d}\, F_{{Q}_l}(\omega)= \gamma(\lambda)$ holds.

Since $S$ is countable, $\Omega'=\cap_{\lambda\in S}\Omega_\lambda$
still has full measure and the convergence statement of the superadditive theorem holds for
all $\lambda \in S $ and $\omega \in \Omega'$.
Now define the distribution function $N(\lambda):= \lim_{S \ni \tilde \lambda \searrow \lambda} \gamma(\tilde \lambda)$.
Thus on the set $\cC$ the functions $\gamma$ and $N$ coincide.

The monotonicity of $ \lambda \mapsto F_{Q_l}(\lambda, \omega)$
and the continuity of  $N$ on the set $\cC$ imply the convergence
\eqref{e-convergence}. To see this, choose a sequence $\lambda_n
\in S,\, \lambda_n \ge \lambda,\, \lim_{n \to \infty} \lambda_n =
\lambda$. Then we have
\[
l^{-d} F_{Q_l} (\lambda,\omega) - N(\lambda)
\le
l^{-d} F_{Q_l} (\lambda_n,\omega) - N(\lambda_n) +N(\lambda_n) -N(\lambda).
\]
For arbitrary $\omega \in \Omega'$ and $\epsilon >0$ we choose first $n$  sufficiently large such that
$N(\lambda_n) -N(\lambda) \le \epsilon/2$ and then $l$ sufficiently large such that
$ l^{-d} F_{Q_l} (\lambda_n,\omega) - N(\lambda_n) \le \epsilon/2$. Thus one sees that
\[
\limsup_{l\to\infty}\, l^{-d} F_{Q_l} (\lambda,\omega) \le N(\lambda).
\]
Similarly one can choose a sequence $\lambda_n \in S,\, \lambda_n
\le \lambda, \,\lim_{n \to \infty} \lambda_n = \lambda$ and then
show that $\liminf_{l\to\infty} l^{-d} F_{Q_l} (\lambda,\omega)
\ge N(\lambda)$. Thus the theorem is proven.
\end{proof}

\bigskip

\noindent
{\textbf{ Acknowledgements:}
Illuminating discussions with D.~Lenz which led to an improvement of the manuscript are gratefully acknowledged.
We thank C.~Schubert for making a draft of the thesis \cite{Schubert-06} available to us.
The second named author was financially supported by the DFG under grant
Ve 253/2-2 within the Emmy-Noether-Programme.}

\medskip

{\textbf{ Note added:} A few days before we finished this manuscript, the paper 
\emph{Anderson Localization for radial tree-like random quantum graphs}
by P.~Hislop and O.~Post, containing results not unrelated to ours, was made available 
at www.arxiv.org.}

\newcommand{\etalchar}[1]{$^{#1}$}
\def\cprime{$'$}


\begin{thebibliography}{DLM{\etalchar{+}}03}

\bibitem[AK81]{AkcogluK-81}
M.~A. Akcoglu and U.~Krengel.
\newblock Ergodic theorems for superadditive processes.
\newblock {\em J. Reine Angew. Math.}, 323:53--67, 1981.

\bibitem[CFKS87]{CyconFKS-87}
H.~L. Cycon, R.~G. Froese, W.~Kirsch, and B.~Simon.
\newblock {\em {Schr\"odinger} Operators with Application to Quantum Mechanics
  and Global Geometry}.
\newblock Text and Monographs in Physics. Springer, Berlin, 1987.

\bibitem[CHK]{CombesHK}
J.-M. Combes, P.~D. Hislop, and F.~Klopp.
\newblock An optimal {Wegner} estimate and its application to the global
  continuity of the integrated density of states for random {Schr\"odinger}
  operators.
\newblock http://arxiv.org/abs/math-ph/0605029.

\bibitem[CHK03]{CombesHK-03}
J.-M. Combes, P.~D. Hislop, and F.~Klopp.
\newblock H\"older continuity of the integrated density of states for some
  random operators at all energies.
\newblock {\em Int. Math. Res. Not.}, (4):179--209, 2003.

\bibitem[CHN01]{CombesHN-01}
J.-M. Combes, P.~D. Hislop, and S.~Nakamura.
\newblock The ${L}^p$-theory of the spectral shift function, the {Wegner}
  estimate, and the integrated density of states for some random
  {Schr\"odinger} operators.
\newblock {\em Commun. Math. Phys.}, 70(218):113--130, 2001.

\bibitem[CL90]{CarmonaL-90}
R.~Carmona and J.~Lacroix.
\newblock {\em Spectral Theory of Random {Schr\"odinger} Operators}.
\newblock Birkh\"auser, Boston, 1990.

\bibitem[DLM{\etalchar{+}}03]{DodziukLMSY-03}
J.~Dodziuk, P.~Linnell, V.~Mathai, T.~Schick, and S.~Yates.
\newblock Approximating {$L^2$}-invariants, and the {Atiyah} conjecture.
\newblock {\em Comm. Pure Appl. Math.}, 56(7):839--873, 2003.

\bibitem[Gie]{Giere-02}
E.~Giere.
\newblock Wegner estimate in one dimension for nonoverlapping single site
  potentials.
\newblock www.math.tu-clausthal.de/$\sim$maeg/arbeiten/publications.html, 2002.

\bibitem[GLV]{GruberLV}
M.~Gruber, D.~Lenz, and I.~Veseli{\'c}.
\newblock Uniform existence of the integrated density of states for random
  {Schr\"odinger} operators on metric graphs over $\mathbb{Z}^d$.
\newblock Preprint.

\bibitem[GT83]{GilbargT-83}
D.~Gilbarg and N.~S. Trudinger.
\newblock {\em Elliptic partial differential equations of second order}.
\newblock Springer, Berlin, 1983.

\bibitem[Kir96]{Kirsch-96}
W.~Kirsch.
\newblock {Wegner} estimates and {Anderson} localization for alloy-type
  potentials.
\newblock {\em Math. Z.}, 221:507--512, 1996.

\bibitem[Klo95]{Klopp-95a}
F.~Klopp.
\newblock Localization for some continuous random {Schr\"odinger} operators.
\newblock {\em Commun. Math. Phys.}, 167:553--569, 1995.

\bibitem[KLS03]{KlassertLS-03}
S.~Klassert, D.~Lenz, and P.~Stollmann.
\newblock Discontinuities of the integrated density of states for random
  operators on {D}elone sets.
\newblock {\em Comm. Math. Phys.}, 241(2-3):235--243, 2003.
\newblock http://arXiv.org/math-ph/0208027.

\bibitem[KM]{KirschM}
W.~Kirsch and B.~Metzger.
\newblock The integrated density of states for random {Schr\"odinger}
  operators.
\newblock http://www.arXiv.org/abs/math-ph/0608066.

\bibitem[KM82]{KirschM-82c}
W.~Kirsch and F.~Martinelli.
\newblock On the density of states of {Schr\"odinger} operators with a random
  potential.
\newblock {\em J. Phys. A: Math. Gen.}, 15:2139--2156, 1982.

\bibitem[Kre85]{Krengel-85}
U.~Krengel.
\newblock {\em Ergodic Theorems}.
\newblock Studies in Mathematics. Walter de Gruyter, 1985.

\bibitem[KS99]{KostrykinS-99b}
V.~Kostrykin and R.~Schrader.
\newblock Kirchhoff's rule for quantum wires.
\newblock {\em J. Phys. A}, 32(4):595--630, 1999.

\bibitem[Kuc04]{Kuchment-04}
P.~Kuchment.
\newblock Quantum graphs. {I}. {S}ome basic structures.
\newblock {\em Waves Random Media}, 14(1):S107--S128, 2004.
\newblock Special section on quantum graphs.

\bibitem[KV02]{KirschV-02b}
W.~Kirsch and I.~Veseli\'c.
\newblock Existence of the density of states for one-dimensional alloy-type
  potentials with small support.
\newblock In {\em Mathematical Results in Quantum Mechanics (Taxco, Mexico,
  2001)}, volume 307 of {\em Contemp. Math.}, pages 171--176. Amer. Math. Soc.,
  Providence, RI, 2002.
\newblock http://arxiv.org/abs/math-ph/0204030.

\bibitem[LPV04]{LenzPV-04}
D.~Lenz, N.~Peyerimhoff, and I.~Veseli\'c.
\newblock Integrated density of states for random metrics on manifolds.
\newblock {\em Proc. London Math. Soc. (3)}, 88(3):733--752, 2004.

\bibitem[LV]{LenzV}
D.~Lenz and I.~Veseli{\'c}.
\newblock Hamiltonians on discrete structures: jumps of the integrated density
  of states.
\newblock in preparation.

\bibitem[PF92]{PasturF-92}
L.~A. Pastur and A.~L. Figotin.
\newblock {\em Spectra of Random and Almost-Periodic Operators}.
\newblock Springer Verlag, Berlin, 1992.

\bibitem[RS78]{ReedS-78}
M.~Reed and B.~Simon.
\newblock {\em Methods of Modern Mathematical Physics {IV}, Analysis of
  Operators}.
\newblock Academic Press, San Diego, 1978.

\bibitem[Sch06]{Schubert-06}
C.~Schubert.
\newblock Laplace-{Operatoren} auf {Quantengraphen}.
\newblock Diplomarbeit, TU Chemnitz, 2006.

\bibitem[Sim95]{Simon-95}
B.~Simon.
\newblock Spectral analysis of rank one perturbations and applications.
\newblock In L.~Rosen J.~Feldman, R.~Froese, editor, {\em Mathematical quantum
  theory. II. Schr\"odinger operators (Vancouver, 1993)}, volume~8 of {\em CRM
  Proc. Lecture Notes Vol. 8}, pages 109--149. Amer. Math. Soc., Providence,
  RI, 1995.

\bibitem[Sto01]{Stollmann-01}
P.~Stollmann.
\newblock {\em Caught by disorder: A Course on Bound States in Random Media},
  volume~20 of {\em Progress in Mathematical Physics}.
\newblock Birkh\"auser, 2001.

\bibitem[Ves96]{Veselic-96}
I.~Veseli\'{c}.
\newblock Lokalisierung bei zuf\"allig gest\"orten periodischen
  {Schr\"odinger\-opera\-toren} in {Dimension Eins}.
\newblock Diplomarbeit, Ruhr-Universit\"at Bochum, 1996.
\newblock available at
  http://www.ruhr-uni-bochum.de/mathphys/ivan/diplomski-www-abstract.htm.

\bibitem[Ves04]{Veselic-04a}
I.~Veseli\'c.
\newblock Integrated density of states and {Wegner} estimates for random
  {Schr\"odinger} operators.
\newblock {\em Contemp. Math.}, 340:98--184, 2004.
\newblock http://arXiv.org/math-ph/0307062.

\bibitem[Ves05]{Veselic-05b}
I.~Veseli\'c.
\newblock Spectral analysis of percolation {Hamiltonians}.
\newblock {\em Math.~Ann.}, 331(4):841--865, 2005.
\newblock http://arXiv.org/math-ph/0405006.

\bibitem[Ves06]{Veselic-06}
I.~Veseli\'c.
\newblock Spectral properties of {Anderson}-percolation {Hamiltonians}.
\newblock {\em Oberwolfach Rep.}, 3(1):545--547, 2006.

\end{thebibliography}
\end{document}